\documentclass[12pt]{article}
\usepackage{amssymb,amsmath,mathrsfs}
\usepackage{hyperref}
\usepackage{graphicx}
\usepackage{color}
\usepackage[OT2,OT1]{fontenc}
\usepackage{upquote}
\usepackage{authblk}
\usepackage[numbers,sort&compress]{natbib}
\usepackage{float}

\begin{document}

\title{\Large\bf A generalized series expansion of the arctangent function based on the enhanced midpoint integration}

\author{
{\small Sanjar M. Abrarov, Rehan Siddiqui, Rajinder K. Jagpal,} \\
{\small and Brendan M. Quine}
}

\date{April 18, 2023}

\maketitle
\maketitle

\begin{abstract}
In this work we derive a generalized series expansion of the acrtangent function by using the enhanced midpoint integration (EMI). Algorithmic implementation of the generalized series expansion utilizes two-step iteration without surd and complex numbers. The computational test we performed reveals that such a generalization improves accuracy in computation of the arctangent function by many orders of the magnitude with increasing integer $M$, associated with subintervals in the EMI formula. The generalized series expansion may be promising for practical applications. It may be particularly useful in practical tasks, where extensive computations with arbitrary precision floating points are needed. The algorithmic implementation of the generalized series expansion of the arctangent function shows a rapid convergence rate in the computation of digits of $\pi$ in the Machin-like formulas.
\vspace{0.25cm}
\vspace{0.25cm}
\end{abstract}

\section{Introduction}

In 2010, Adegoke and Layeni published an interesting relation for derivatives of the arctangent function \cite{Adegoke2010}
\begin{equation}\label{eq_1}
\frac{d^n}{dx^n}\arctan(x) = \frac{(-1)^{n - 1}(n - 1)!}{(1 + x^2)^{n/2}}\sin\left(n\arcsin\left(\frac{1}{\sqrt{1 + {x^2}}}\right)\right), \,\, n \in \mathbb N^+.
\end{equation}
Using this relation, they discovered a series expansion:
\begin{equation}\label{eq_2}
\arctan(x) = \sum\limits_{n = 1}^\infty\frac{1}{n}\left(\frac{x^2}{1 + x^2}\right)^{n/2}\sin\left(n\arcsin\left(\frac{1}{\sqrt{1 + x^2}}\right)\right).
\end{equation}

Equations \eqref{eq_1} and \eqref{eq_2} have some restrictions. Specifically, when $n$ is even, equation \eqref{eq_1} remains valid only at $x \in [0,\infty)$, while equation \eqref{eq_2} is valid only at $x \in [0,\infty)$ for $\forall n$.

To resolve this problem, Lampret applied the signum function
$$
{\rm sgn}(x) = \left\{
\begin{aligned}
-1, \qquad x < 0, \\
1, \qquad x \ge 0, \\ 
\end{aligned}
\right.
$$
and proved that for complete coverage $x \in(-\infty,\infty)$, the equations \eqref{eq_1} and \eqref{eq_2} can be modified as \cite{Lampret2011}
\small
\begin{equation}\label{eq_3}
\frac{d^n}{dx^n}\arctan(x) = {\rm sgn}(-x)^{n - 1}\frac{(n - 1)!}{(1 + x^2)^{n/2}}\sin\left(n\arcsin\left(\frac{1}{\sqrt{1 + {x^2}}}\right)\right)
\end{equation}
\normalsize
and
\small
\begin{equation}\label{eq_4}
\arctan(x) = {\rm sgn}(x)\sum\limits_{n = 1}^\infty\frac{1}{n}\left(\frac{x^2}{1 + x^2}\right)^{n/2}\sin\left(n\arcsin\left(\frac{1}{\sqrt{1 + x^2}}\right)\right),
\end{equation}
\normalsize
respectively.

Equations \eqref{eq_3} and \eqref{eq_4} represent a theoretical interest. In particular, Lampret noticed that from equation \eqref{eq_3}, it follows that \cite{Lampret2011}
\begin{equation}\label{eq_5}
{\rm sgn}(0)^{n - 1} \cdot (n - 1)! \cdot \sin\left(n\frac{\pi}{2}\right) =
\left\{\begin{aligned}
&(-1)^{(n - 1)/2}(n - 1)!, &&\,\, n\,\,{\rm odd}, \\
&0, &&\,\, n\,\,{\rm even}. \\ 
\end{aligned}\right.
\end{equation}

Comparing the following relation (see \cite{Abrarov2016a} for detailed derivation procedure by induction)
\begin{equation}\label{eq_6}
\frac{d^n}{dx^n}\arctan(x) = \frac{(-1)^n(n - 1)!}{2i}\left(\frac{1}{(x + i)^n} - \frac{1}{(x - i)^n}\right)
\end{equation}
with equation \eqref{eq_3}, we can find the following identity:
\begin{equation}\label{eq_7}
\begin{aligned}
{\rm sgn}(-x)^{n - 1}\frac{(n - 1)!}{(1 + x^2)^{n/2}}\sin&\left(n\arcsin\left(\frac{1}{\sqrt{1 + x^2}}\right)\right) = \\
&\frac{(-1)^n(n-1)!}{2i}\left(\frac{1}{(x + i)^n} - \frac{1}{(x - i)^n}\right).
\end{aligned}
\end{equation}

It is not difficult to see that the relation \eqref{eq_5} immediately follows from the identity \eqref{eq_7}. Therefore, relation \eqref{eq_5} is just a specific case of the identity \eqref{eq_7} occurring at $x = 0$.

Identity \eqref{eq_7} can be rewritten in form
\small
\[
\sin\left(n\arcsin\left(\frac{1}{\sqrt{1 + x^2}}\right)\right) = -{\rm sgn}(x)^{n-1}\frac{(1 + x^2)^{n/2}}{2i}\left(\frac{1}{(x + i)^n} - \frac{1}{(x - i)^n}\right).
\]
\normalsize
Therefore, from equation \eqref{eq_4}, it follows that
\small
\[
\begin{aligned}
&\arctan(x) = \\
& \enspace\quad{\rm sgn}(x)\sum_{n = 1}^\infty\frac{1}{n}\left(\frac{x^2}{1 + x^2}\right)^{n/2}\left[-{\rm sgn}(x)^{n-1}\frac{(1 + x^2)^{n/2}}{2i}\left(\frac{1}{(x + i)^n} - \frac{1}{(x - i)^n}\right)\right]
\end{aligned}
\]
\normalsize
or
\[
\arctan(x) = \frac{i}{2}\sum_{n = 1}^\infty\frac{x^n}{n}\left(\frac{1}{(x + i)^n} - \frac{1}{(x - i)^n}\right).
\]
As we can see, this series expansion of the arctangent function is just a reformulation of equation \eqref{eq_4} and, since the relation \eqref{eq_6} can be rearranged in form
\[
\frac{1}{(x + i)^n} - \frac{1}{(x - i)^n} = \frac{2i}{(-1)^n(n - 1)!}\frac{d^n}{dx^n}\arctan(x),
\]
we can express the arctangent function in terms of its derivatives as given by the following equation:
$$
\arctan(x) = \sum\limits_{n = 1}^\infty\frac{(-1)^{n - 1}x^n}{n!}\frac{d^n}{d x^n}\arctan(x)
$$
that leads to
$$
\sum\limits_{n = 0}^\infty\frac{(-1)^{n - 1}x^n}{n!}\frac{d^n}{d x^n}\arctan(x) = 0.
$$

In our previous publication \cite{Abrarov2017a}, using the identity \eqref{eq_6}, we have derived the following series expansion of the arctangent function
\[
\begin{aligned}
&\arctan(x) =  \\
&\hspace{0.75cm} -2\sum\limits_{m = 1}^\infty\sum\limits_{n = 1}^{2m - 1}\frac{(-1)^n}{(2m-1)(1 + x^2/4)^{2m - 1}}\left(\frac{x}{2}\right)^{2(2m - n) - 1}{2m - 1 \choose 2n - 1},
\end{aligned}
\]
from which, at $x = 1$, we get a formula for $\pi$ expressed in terms of the binomial coefficients
\[
\frac{\pi}{4} =  -2\sum\limits_{m = 1}^\infty\sum\limits_{n = 1}^{2m - 1}\frac{(-1)^n}{(2m - 1)(1 + 1/4)^{2m - 1}2^{2(2m - n) - 1}}{2m - 1 \choose 2n - 1} 
\]
or
\[
\frac{\pi}{16} = \sum\limits_{m = 1}^\infty\sum\limits_{n = 1}^{2m - 1}\frac{(-4)^{n - 1}}{(2m - 1)5^{2m - 1}}{2m - 1 \choose 2n - 1}.
\]

Later, using the same identity \eqref{eq_6}, we have also derived the following series expansion (see \cite{Abrarov2017b,Abrarov2022})
\begin{equation}\label{eq_8}
\arctan(x) = 2\sum\limits_{n = 1}^\infty\frac{1}{2n - 1}\frac{g_n(x)}{g_n^2(x) + h_n^2(x)},
\end{equation}
where the expansion coefficients are computed by two-step iteration:
\[
\begin{gathered}
g_n(x) = g_{n - 1}(x)(1 - 4/x^2) + 4h_{n - 1}(x)/x, \\
h_n(x) = h_{n - 1}(x)(1 - 4/x^2) - 4g_{n - 1}(x)/x, 
\end{gathered}
\]
such that
\[
\begin{gathered}
g_1(x) = 2/x, \\
h_1(x) = 1. 
\end{gathered}
\]
The series expansion \eqref{eq_8} requires no surd or complex numbers in computation and it is rapid in convergence.

Many new identities and series expansions related to the arctangent function have been reported \cite{Milgram2006,Sofo2012,Zhang2013,Qi2015,Pilato2017,Qiao2018,Benammar2019,Sofo2020,Kusaka2022}. This shows that the discovery of new equations related to the arctangent function as well as their applications remain a very interesting topic.

As further development, in this work we derive a generalized series expansion of the arctangent function. Such an approach may be used to improve further convergence in computation of the arctangent function. Due to rapid convergence without surd and complex numbers in computation, the generalized series expansion may be promising for applications with arbitrary precision floating points \cite{Brent1978,Brent2006,Dinechin2006,Bailey2012,Bailey2015,Kneusel2017,Vestermark2022,Johansson2022}. Furthermore, it may also be promising in computing digits of $\pi$ by using the Machin-like formulas \cite{Abrarov2022,Berggren2004,Calcut2009,Nimbran2010,Agarwal2013,Wetherfield,Abrarov2017b,WolframCloud,Gasull2023}. To the best of our knowledge, the generalized series expansion of the arctangent function is new and has never been reported.

\section{Derivation}

Change of the variable $x \to xt$ in the equation \eqref{eq_6} results in
\[
\frac{\partial^n}{\partial t^n}\arctan(xt) = \frac{(-1)^n(n - 1)!x^n}{2i}\left(\frac{1}{(xt + i)^n} - \frac{1}{(xt - i)^n}\right)
\]
or
\begin{equation}\label{eq_9}
\frac{\partial^n}{\partial{t^n}}\frac{x}{1 + x^2t^2} = \frac{(-1)^{n + 1}n!x^{n + 1}}{2i}\left(\frac{1}{(xt + i)^{n + 1}} - \frac{1}{(xt - i)^{n + 1}}\right)
\end{equation}
since
\[
\frac{\partial^n}{\partial{t^n}}\arctan(xt) = \frac{\partial^{n - 1}}{\partial{t^{n - 1}}}\left(\frac{x}{1 + x^2t^2}\right).
\]

The Enhanced Midpoint Integration (EMI) formula is given by (see \cite{Abrarov2016b} for derivation and \cite{Abrarov2018} for application)
\begin{equation}\label{eq_10}
\int\limits_0^1 f(t)\,dt = \sum\limits_{m = 1}^M\sum\limits_{n = 0}^\infty\frac{(-1)^n + 1}{(2M)^{n + 1}(n + 1)!}{\left.\frac{d^n}{d{t^n}}f(t)\right|}_{t = \frac{m - 1/2}{M}},
\end{equation}
where integer $M$ is associated with subintervals of integration. It is interesting to note that if the upper summation bound associated with variable $n$ is an integer $N \ge 0$, then we can also use
$$
\int\limits_0^1 f(t)\,dt = \lim\limits_{M\to\infty}\sum\limits_{m = 1}^M\sum\limits_{n = 0}^N\frac{(-1)^n + 1}{(2M)^{n + 1}(n + 1)!}{{\left.\frac{d^n}{d{t^n}}f(t)\right|}_{t = \frac{m - 1/2}{M}}}.
$$

It is easy to show that, excluding all zero terms occurring at odd values of the variable $n$, equation \eqref{eq_10} can be rewritten in a more convenient form:
\begin{equation}\label{eq_11}
\int\limits_0^1 f(t)\,dt = 2\sum\limits_{m = 1}^M\sum\limits_{n = 0}^\infty\frac{1}{(2M)^{2n + 1}(2n + 1)!}{{\left.\frac{d^{2n}}{d{t^{2n}}}f(t)\right|}_{t = \frac{m - 1/2}{M}}}.
\end{equation}

Equation \eqref{eq_11} requires even derivatives of the integrand at the points $t = (m - 1/2)/M$, where $m = 1,2,3,\,\, \ldots \, M$. Taking these derivatives manually is extremely tedious. However, with help of the Computer Algebra System (CAS) such as Mathematica, Matlab or Maple supporting symbolic programming, the application of the equation \eqref{eq_11} may be very efficient in numerical integration. Specifically, such an approach may be especially useful for the numerical integration of the highly oscillating functions.

The EMI formula \eqref{eq_11} can be used for numerical integration within any interval $t \in (a,b)$, since the following transformation formula
\[
\int\limits_a^b f(t)\,dt = (b - a)\int\limits_0^1 f((b - a)t + a)\,dt
\]
can be applied to recast the integration interval within $t \in (0,1)$. The interested readers can download the MATLAB code based on the integration formula \eqref{eq_11} on the MATLAB Central website \cite{MatlabCentral} (file ID \#: 71037). The synopsis and brief instruction on how to use the MATLAB code for numerical integration is also provided in the supplementary {\it readme.pdf} file.

If an integrand represents a function of two variables $f\left( {x,t} \right)$, then the integration formula \eqref{eq_11} reads as
\begin{equation}\label{eq_12}
\int\limits_0^1 f(x,t)\,dt = 2\sum\limits_{m = 1}^M\sum\limits_{n = 0}^\infty\frac{1}{(2M)^{2n + 1}(2n + 1)!}{{\left.\frac{\partial^{2n}}{\partial{t^{2n}}}f(x,t)\right|}_{t = \frac{m - 1/2}{M}}}.
\end{equation}

The arctangent function can be given as an integral
\begin{equation}\label{eq_13}
\arctan(x) = \int\limits_0^1\frac{x}{1 + x^2t^2}\,dt.
\end{equation}
Consequently, substituting the integrand from equation \eqref{eq_13}
$$
\frac{x}{1 + x^2t^2} = \frac{x}{2}\left(\frac{1}{1 + ixt} + \frac{1}{1 - ixt}\right)
$$
into equation \eqref{eq_12} and using equation \eqref{eq_9} for differentiation, we can find that
\small
\begin{equation}\label{eq_14}
\begin{aligned}
&\arctan(x) = \\
&\enspace\quad i\sum\limits_{m = 1}^M\sum\limits_{n = 0}^\infty\frac{x^{2n + 1}}{(2M)^{2n + 1}(2n + 1)}\left[\frac{1}{(x\frac{{m - 1/2}}{M} + i)^{2n + 1}} - \frac{1}{(x\frac{m - 1/2}{M} - i)^{2n + 1}}\right].
\end{aligned}
\end{equation}
\normalsize

Series expansion \eqref{eq_14} is rapid in convergence. However, it requires algebraic manipulations with complex numbers. Therefore, it is very desirable to exclude them. This can be achieved by induction 
based on two-step iteration
$$
\begin{gathered}
\alpha_n(x,t) = \alpha_{n - 1}(x,t)(1 - 1/(xt)^2) + 2\beta_{n - 1}(x,t)/(xt), \\
{\beta _n}\left( {x,t} \right) = {\beta _{n - 1}}\left( {x,t} \right)\left( {1 - 1/{{\left( {xt} \right)}^2}} \right) - 2{\alpha _{n - 1}}\left( {x,t} \right)/\left( {xt} \right), \\ 
\end{gathered}
$$
and
$$
\begin{gathered}
\alpha_1(x,t) = 1/(xt), \\
\beta_1(x,t) = 1, \\ 
\end{gathered}
$$
that transforms equation \eqref{eq_14} into the following series expansion:
\small
\begin{equation}\label{eq_15}
\arctan(x) = 2\sum\limits_{m = 1}^M\sum\limits_{n = 1}^\infty\frac{1}{(2n - 1)(2m - 1)^{2n - 1}}\,\frac{\alpha_n(x,\gamma_{m,M})}{\alpha_n^2(x,\gamma_{m,M}) + \beta_n^2(x,\gamma_{m,M})},
\end{equation}
\normalsize
where the argument is
$$
\gamma_{m,M} = \frac{m - 1/2}{M}.
$$

Equation \eqref{eq_15} is a generalization of the equation \eqref{eq_8}. Consistency between these two equations can be observed by taking $M = 1$. In particular, substitution $M = 1$ into series expansion \eqref{eq_15} of the arctangent function implies that $t = \gamma_{1,1} = 1/2$. Consequently, from equation \eqref{eq_15} we obtain equation \eqref{eq_8}, where the expansion coefficients are
$$
\begin{gathered}
g_1(x) = \alpha_1(x,\gamma_{1,1}), \\
h_1(x) = \beta_1(x,\gamma_{1,1}), \\
g_n(x) = \alpha_n(x,\gamma_{1,1}), \\
h_n(x) = \beta_n(x,\gamma_{1,1}). \\ 
\end{gathered}
$$
\begin{figure}[H]
\begin{center}
\includegraphics[width=24pc]{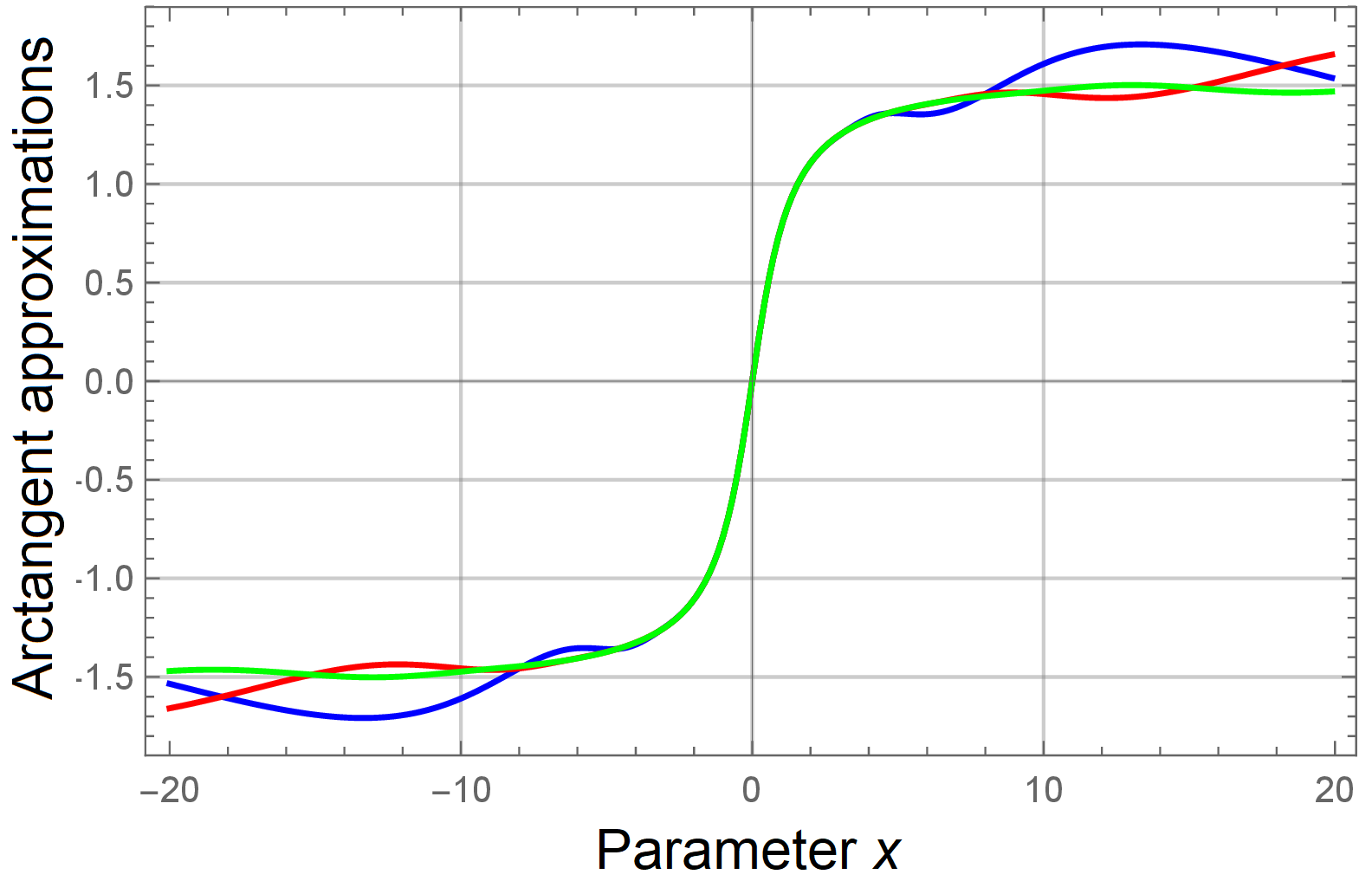}\hspace{2pc}%
\begin{minipage}[b]{28pc}
\vspace{0.3cm}
{\sffamily {\bf{Fig. 1.}} Arctangent approximations computed by using series expansion \eqref{eq_15} truncated at $n_{\rm max} = 10$. Blue, red and green curves correspond to $M$ taken to be $1$, $2$ and $3$, respectively.}
\end{minipage}
\end{center}
\end{figure}

The following is a Mathematica code that generates graphs shown in Fig.~1 (this code can be copy-pasted directly to the Mathematica notebook):

\pagestyle{empty}
\footnotesize
\begin{verbatim}
Clear[atan,\[Gamma],\[Alpha],\[Beta]];

(* Equation (15) *)
atan[x_,nMax_,M_] := 2*Sum[(1/(2*n - 1))*
  (\[Alpha][x,\[Gamma][m,M],n]/((2*m - 1)^(2*n - 1)*
    (\[Alpha][x,\[Gamma][m, M],n]^2 + \[Beta][x,
      \[Gamma][m,M],n]^2))),{m,1,M},{n,1,nMax}];

(* Argument gamma *)
\[Gamma][m_,M_] := \[Gamma][m,M] = N[(m - 1/2)/M,1000];

(* Expansion coefficients *)
\[Alpha][x_,t_,1] := \[Alpha][x,t,1] = 1/(x*t);
\[Beta][x_,t_,1] := \[Beta][x,t,1] = 1;

\[Alpha][x_,t_,n_] := \[Alpha][x,t,n] =
  \[Alpha][x,t,n - 1]*(1 - 1/(x*t)^2) +
    2*(\[Beta][x,t,n - 1]/(x*t));

\[Beta][x_,t_,n_] := \[Beta][x,t,n] = 
  \[Beta][x,t,n - 1]*(1 - 1/(x*t)^2) -
    2*(\[Alpha][x,t,n - 1]/(x*t));

(* Computing data points *)
tabs := {Table[{x,atan[x,10,1]},{x,-20,20,Pi/20}],
  Table[{x,atan[x,10,2]},{x,-20,20,Pi/20}],
    Table[{x,atan[x,10,3]},{x,-20,20,Pi/20}]};

Print["Computing, please wait..."];

(* Plotting graphs *)
ListPlot[tabs,Joined->True,FrameLabel->{"Parameter x", 
  "Arctangent approximations"},PlotStyle->{Blue,Red,Green}, 
    Frame->True,GridLines->Automatic]
\end{verbatim}
\normalsize
\pagestyle{plain}

The graphs in Fig. 1 are generated by using series expansion \eqref{eq_15} truncated at $n_{\rm max} = 10$. Blue, red and green curves correspond to integer $M$ taken to be $1$, $2$ and $3$, respectively.

Consider Fig. 2 showing approximation curves of the arctangent function $\arctan(x)$ by using equations \eqref{eq_15}, \eqref{eq_16} and \eqref{eq_17} truncated at $n_{\rm max} = 10$. The blue curve corresponding to the Maclaurin series expansion
\begin{equation}\label{eq_16}
\arctan(x) = \sum\limits_{n = 0}^\infty \left.{\frac{x^n}{n!}\frac{d^n}{dt^n}\arctan(t)}\right|_{t = 0} = \sum\limits_{n = 0}^\infty\frac{(-1)^n x^{2n + 1}}{2n + 1}, \qquad|x| \le 1,
\end{equation}
diverges beyond $-1$ and $1$ due to finite radius of convergence. Although one can resolve this issue by using an elementary relation
$$
\arctan \left( x \right) = \left\{
\begin{aligned}
&\frac{\pi}{2} - \arctan\left(\frac{1}{x}\right), &&\quad x > 0, \\
&0, && \quad x = 0, \\
&-\frac{\pi}{2} - \arctan\left(\frac{1}{x}\right), &&\quad x < 0, \\ 
\end{aligned}
\right.
$$
our objective is just to visualize the convergence. The red curve shows the Euler series expansion \cite{Vestermark2022,Castellanos1988,Chien-Lih2005}:
\begin{equation}\label{eq_17}
\arctan \left( x \right) = \sum\limits_{n = 0}^\infty  {\frac{{{2^{2n}}{{\left( {n!} \right)}^2}}}{{\left( {2n + 1} \right)!}}\frac{{{x^{2n + 1}}}}{{{{\left( {1 + {x^2}} \right)}^{n + 1}}}}.}
\end{equation}
The green curve illustrates the series expansion \eqref{eq_15} at $M = 1$. The black dashed curve depicts the original arctangent function for comparison. As we can see from Fig. 2, even at smallest $M$ the series expansion \eqref{eq_15} provides more rapid convergence as compared to the Euler series expansion \eqref{eq_17}.
\begin{figure}[ht]
\begin{center}
\includegraphics[width=24pc]{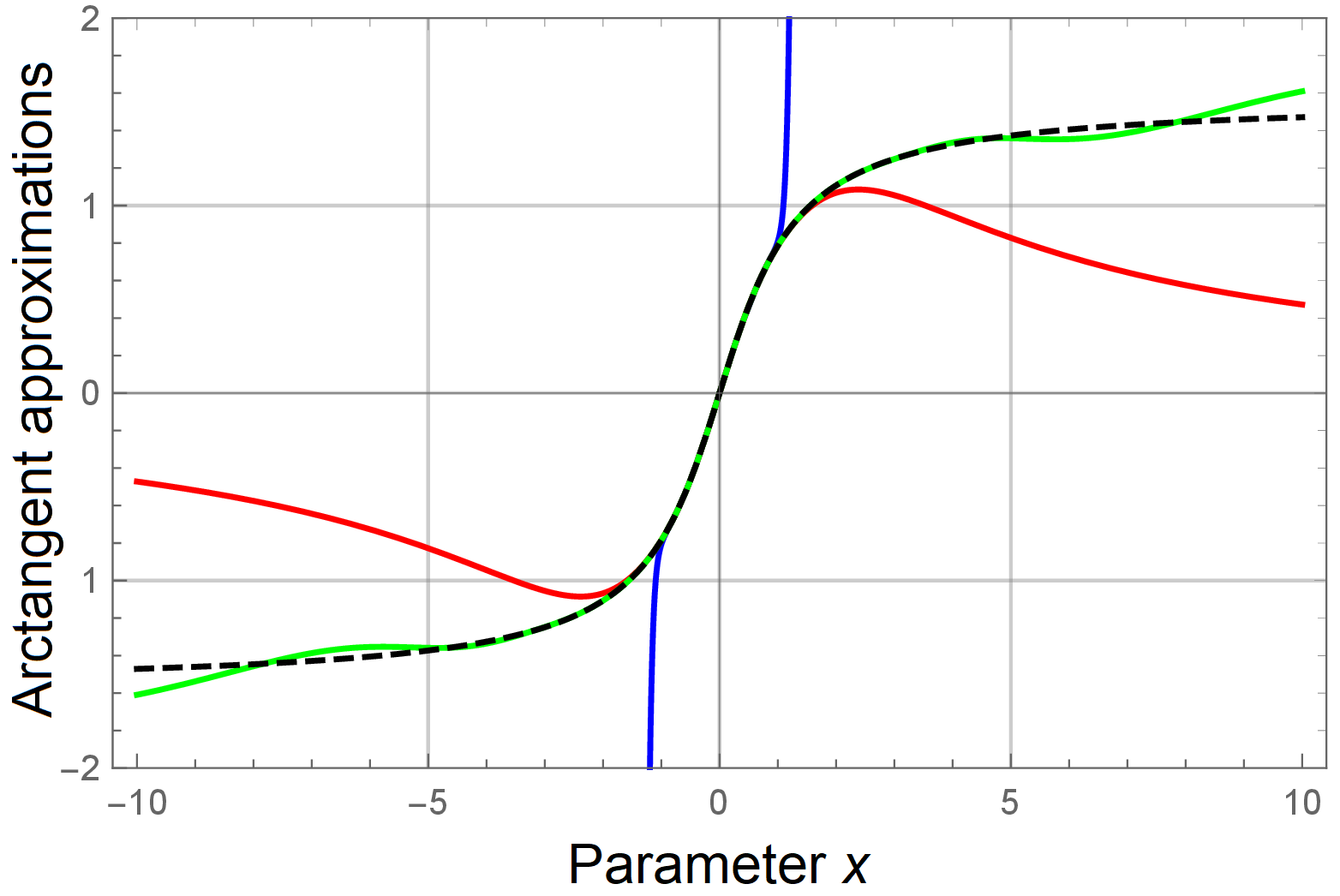}\hspace{2pc}%
\begin{minipage}[b]{28pc}
\vspace{0.3cm}
{\sffamily {\bf{Fig. 2.}} Arctangent approximations computed by series expansions \eqref{eq_16} (blue), \eqref{eq_17} (red) and \eqref{eq_15} (green) truncated at $n_{\rm max} = 10$. Integer $M$ in equation \eqref{eq_15} is taken to be equal to $1$. The dashed black curve shows the original artangent function.}
\end{minipage}
\end{center}
\end{figure}

Figure 3 shows the logarithms of absolute difference $\log_{10}\Delta$ between the arctangent function and its approximations provided by equations \eqref{eq_15}, \eqref{eq_16} and \eqref{eq_17}. All curves are also computed with truncating integer $n_{\rm max} = 10$ in all these equations. The blue and red curves correspond to equations \eqref{eq_16} and \eqref{eq_17} while the green, brown, gray, magenta and black curves correspond to equation \eqref{eq_15} when $M$ is equal to $1$, $2$, $3$, $4$ and $5$, respectively. As we can see from this figure, the increase of the integer $M$ leads to a rapid decrease of the absolute difference $\Delta$ by many orders of the magnitude. These results indicate that the series expansion \eqref{eq_15} provides increasing convergence with increasing $M$.

\begin{figure}[ht]
\begin{center}
\includegraphics[width=24pc]{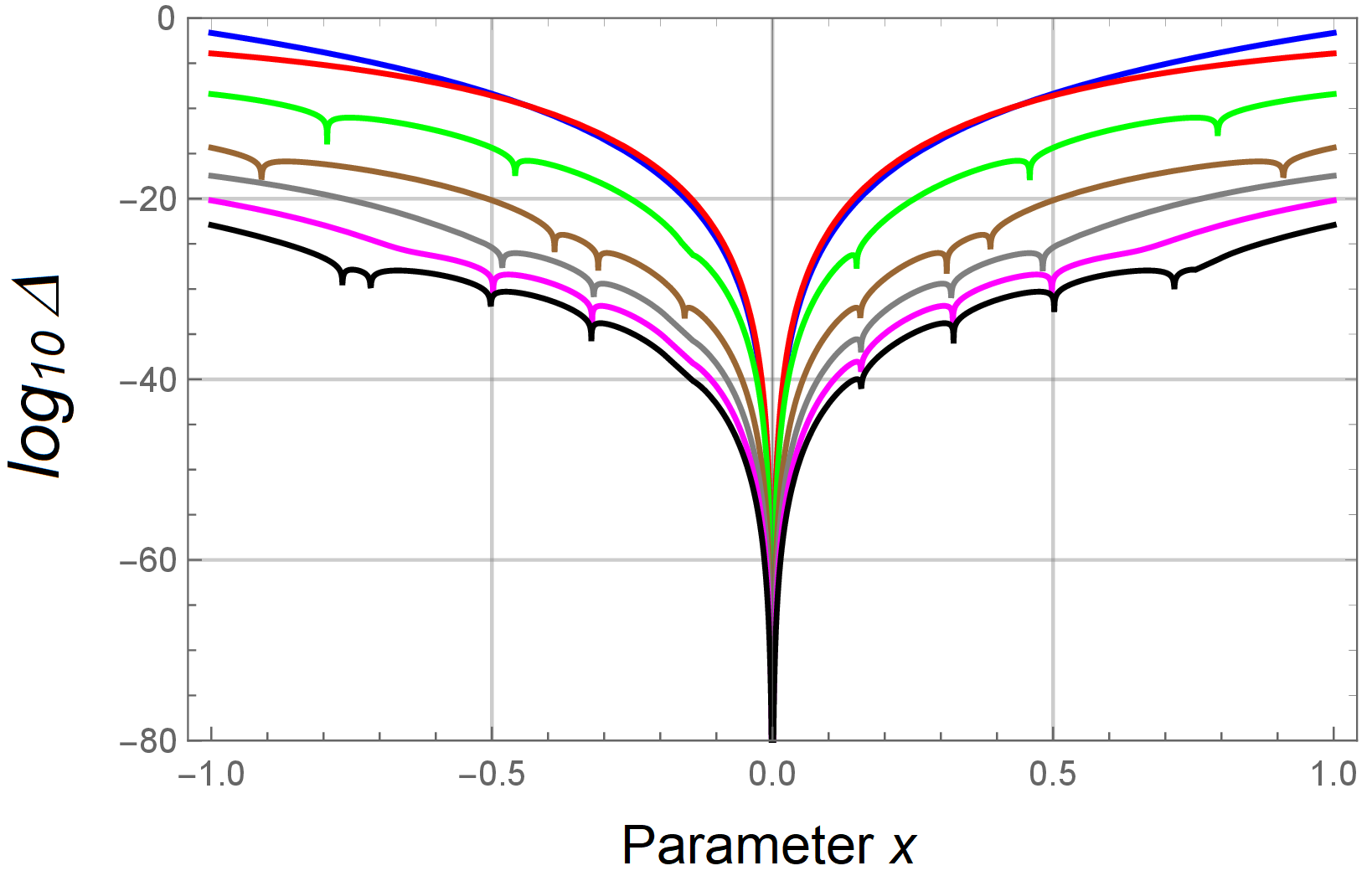}\hspace{2pc}%
\begin{minipage}[b]{28pc}
\vspace{0.3cm}
{\sffamily {\bf{Fig. 3.}} Logarithms of absolute difference $log_{10}\Delta$ between original arctangent function and series expansions \eqref{eq_16} (blue), \eqref{eq_17} (red) and \eqref{eq_15} (green-to-black) truncated at $n_{\rm max} = 10$. Integer $M$ in the series expansion \eqref{eq_15} is taken to be $1$ (green), $2$ (brown), $3$ (gray), $4$ (magenta) and $5$ (black).}
\end{minipage}
\end{center}
\end{figure}

\section{Applications}

There are two potentially possible applications of the series expansion \eqref{eq_15}. In particular, the series expansion \eqref{eq_15} can be implemented as a function file for computation of arctangent for the libraries with arbitrary precision floating points that are required in many fields of Mathematical/Computational Physics and Astronomy \cite{Brent1978,Brent2006,Dinechin2006,Bailey2012,Bailey2015,Kneusel2017,Vestermark2022,Johansson2022}. According to \cite{Brent1978,Brent2006,Vestermark2022,Johansson2022} the Maclaurin and Euler series expansions \eqref{eq_16} and \eqref{eq_17} are used for the arctangent function to provide arbitrary precision in the programing languages.

The Maclaurin series expansion \eqref{eq_16} diverges at the points of argument $1$ and $-1$. This deteriorates accuracy as argument $x$ approaches $1$ by absolute value. Consequently, its algorithmic implementation requires additional algebraic manipulations to overcome this issue for computations with enhanced precision \cite{Brent1978}.

As equation \eqref{eq_15} is more rapid in convergence than equations \eqref{eq_16} and \eqref{eq_17}, its application may also be efficient to reduce the run-time in many tasks, where extensive computations (without surd or complex numbers) with largely extended precision in floating point numbers are necessary. Moreover, equation \eqref{eq_15} may also provide additional flexibility for users, who can choose the parameter $M$ according to their specific requirements for high-accuracy computation.

Another application, where the series expansion \eqref{eq_15} can also find its practical implementation, is a computation of digits of the constant $\pi$ by using the Machin-like formulas \cite{Abrarov2022,Berggren2004,Calcut2009,Nimbran2010,Agarwal2013,Wetherfield,Abrarov2017b,WolframCloud,Gasull2023}.

Since computation of any irrational numbers is itself a big challenge, a rapid convergence of the arctangent terms in the Machin-like formulas without undesirable surd numbers can provide a significant advantage. Our empirical results show that even using already known Machin-like formulas with sufficiently large integers in actangent arguments, the expansion series \eqref{eq_15} at any $M \ge 1$ can provide more than $17$ digits of $\pi$ at each increment by $1$ of the variable $n$. It is interesting to note that this convergence rate is faster than that provided by Chudnovsky formula generating $14$ to $16$ digits of $\pi$ per increment \cite{Agarwal2013,Berggren2004}. Currently, Chudnovsky formula remains most efficient for computing digits of $\pi$ due to its rapid convergence and other advantages in algorithmic implementation. Historically, however, there were several records that appeared due to the application of the Machin-like formulas in computing $\pi$ and, in 2002, an algorithm, developed by Kanada on the basis of a self-checking pair of the Machin-like formulas, beat the record, providing more than a trillion digits of $\pi$ for the first time \cite{Calcut2009,Agarwal2013}. Therefore, the discovery of new Machin-like formulas and rapidly convergent series expansions of the arctangent function may be promising for computing digits of the constant~$\pi$.

Consider the following example. Previously, we developed a method and generated a two-term Machin-like formula for $\pi$ by using two-step iteration process \cite{Abrarov2017b}
\begin{equation}\label{eq_18}
\frac{\pi}{4} = 2^{26}\arctan\left(\frac{1}{85445659}\right) - \arctan\left(\frac{\overbrace{9732933578 \ldots 4975692799}^{522,185,807\,\rm{digits}}}{\underbrace{2368557598 \ldots 9903554561}_{522,185,816\,\rm{digits}}}\right).
\end{equation}
Recently the same identity has been generated by independent researchers, who developed a different method for generating the two-term Machin-like formulas for $\pi$ \cite{Gasull2023}. The method we developed in \cite{Abrarov2017b} for generating two-terms Machin-like formulas for $\pi$ is a complete alternative to the method described in \cite{Gasull2023} for the formulas of kind
\[
\frac{\pi}{4} = 2^{k - 1}\arctan\left(\frac{a}{b}\right) + \arctan\left(\frac{c}{d}\right), \qquad a,b,c,d \in \mathbb{Z}.
\]
Specifically, all Machin-like formulas for $\pi$, shown in the Tables 1 and 2 from \cite{Gasull2023}, can also be generated by using two-step iterative method that we proposed in our publication \cite{Abrarov2017b}. The two-term Machin-like formula \eqref{eq_18} for $\pi$ that is listed in the Table 2 in the recent publication \cite{Gasull2023}, is just a specific case occurring at $k - 1 = 26$, $a = 1$ and $b = 85445659$.

Due to large number of digits in the numerator and denominator in the second arctangent term of  equation \eqref{eq_18}, we cannot perform computation to observe the corresponding convergence rate without a powerful computer. However, this problem can be handled on a typical laptop or desktop computer with the help of the following identity (see derivation in \cite{Abrarov2017b})
\begin{equation}\label{eq_19}
\frac{\pi}{4} = 2^{k - 1}\arctan\left(\frac{1}{\gamma}\right) + \arctan\left(\frac{1 - \sin\left(2^{k - 1}\arctan\left(\frac{2\gamma}{\gamma^2 - 1}\right)\right)}{\cos\left(2^{k - 1}\arctan\left(\frac{2\gamma}{\gamma^2 - 1}\right)\right)}\right),
\end{equation}
where the constant $\gamma$ may be chosen such that
\begin{equation}\label{eq_20}
\frac{2^{k - 1}}{\gamma}\approx\frac{\pi}{4}.
\end{equation}

Since \cite{Abrarov2018}
\[
\frac{\pi}{4} = 2^{k - 1}\arctan\left(\frac{\sqrt{2 - a_{k - 1}}}{a_k}\right), \qquad k\ge 1,
\]
where $a_0 = 0$ and $a_{k + 1} = \sqrt{2 + a_k}$, it is convenient to choose $\gamma$ for equation  as an integer
$$
\gamma = \left\lfloor\frac{a_k}{\sqrt{2 - a_{k - 1}}}\right\rfloor
$$
to satisfy the condition \eqref{eq_20}. Consequently, at $k = 27$ we obtain the value $85445659$ for $\gamma$ that is present in equation \eqref{eq_18} and the difference
\[
\frac{2^{26}}{85445659} - \frac{\pi}{4} = 4.10922\ldots\times 10^{-9}
\]
is small. More generally, the value of $\gamma$ may be taken as a ratio rather than an integer. For example, we can use \cite{Abrarov2022}
$$
\gamma = \left\lfloor 10^m\frac{a_k}{\sqrt{2 - a_{k - 1}}}\right\rfloor 10^{-m}
$$
or
$$
\gamma = \left\lceil 10^m\frac{a_k}{\sqrt{2 - a_{k - 1}}}\right\rceil 10^{-m},
$$
where $m\in\mathbb{N}^+$.

The identity \eqref{eq_19} implies that if the argument of the first arctangent function is known, then the argument of the second arctangent function can be found. Thus, if we take $\gamma = 85445659$, then the argument of the second arctangent function in equation \eqref{eq_18} can be calculated accordingly as
\[
\frac{1 - \sin\left(2^{k - 1}\arctan\left(\frac{2\gamma}{\gamma^2 - 1}\right)\right)}{\cos\left(2^{k - 1}\arctan\left(\frac{2\gamma}{\gamma^2 - 1}\right)\right)} = -4.10922393614549022091\ldots \times 10^{-9}.
\]
Once we obtain the argument of the second arctangent with, say, up to $1000$ correct decimal digits, we can see that substitution of both arctangent arguments into equation \eqref{eq_15} at $M = 1$ gives $15$ to $17$ correct digits of $\pi$ per each increment of $n$. This example demonstrates a rapid convergence rate of the generalized series expansion \eqref{eq_15}. Therefore, its algorithmic implementation may be promising for rapid and highly accurate computation.

\section{Conclusion}

We derived a generalized series expansion \eqref{eq_15} of the arctanget function by using the EMI formula \eqref{eq_12}. Algorithmic implementation of equation \eqref{eq_15} is based on two-step iteration without surd and complex numbers. The computational test we performed reveals that such a generalization significantly improves convergence in computation of the arctangent function with increasing integer $M$. The generalized series expansion \eqref{eq_15} may be promising in practical applications; it may be used for extensive computations with arbitrary precision floating points and its algorithmic implementation shows a high convergence rate in computation of digits of $\pi$ in the Machin-like formulas.

\section*{Acknowledgment}

\addcontentsline{toc}{section}{Acknowledgment}
This work is supported by National Research Council Canada, Thoth Technology Inc., York University and Epic College of Technology.

\bigskip

\end{document}